\documentclass[a4paper,11pt,reqno]{amsart}

\usepackage{amsmath}
\usepackage{amssymb}
\usepackage{hyperref}
\usepackage{todonotes}
\usepackage{algorithm}
\usepackage[noend]{algorithmic}
\usepackage{enumitem}
\usepackage{comment}

\theoremstyle{plain}
\newtheorem{theorem}{Theorem}[section]
\theoremstyle{definition}
\newtheorem{df}[theorem]{Definition}
\newtheorem{example}[theorem]{Example}%
\newtheorem{notation}[theorem]{Notation}%
\newtheorem{rem}[theorem]{Remark}
\theoremstyle{plain}

\theoremstyle{plain}
\newtheorem{prop}[theorem]{Proposition}
\newtheorem{lem}[theorem]{Lemma}

\begin{document}

\title{A smoothness test for higher codimensions}

\author[J.~B\"ohm]{Janko~B\"ohm}
\address{Department of Mathematics\\
University of Kaiserslautern\\
Erwin-Schr\"odinger-Str.\\
67663 Kaiserslautern\\
Germany}
\email{boehm@mathematik.uni-kl.de}

\author[A. Fr\"uhbis-Kr\"uger]{Anne Fr\"uhbis-Kr\"uger}
\address{Institut f\"ur algebraische Geometrie\\
Leibniz Universit\"at Hannover\\
Welfengarten 1\\
30167 Hannover\\
Germany}

\email{anne@math.uni-hannover.de}

\keywords{%
singularities, algebraic geometry, algorithmic smoothness test, Hironaka resolution of singularities, unprojection%
}

\subjclass[2010]{14B05 (Primary), 68W10, 13P10, 32S05 (Secondary)}

\maketitle

\begin{abstract}
Based on an idea in Hironaka's proof of resolution of singularities, we 
present an algorithmic smoothness test for algebraic varieties. The test is 
inherently parallel and does not involve the calculation of codimension-sized 
minors of the Jacobian matrix of the variety. We also describe a hybrid method 
which combines the new method with the Jacobian criterion, thus making use of the strengths of both approaches. We have implemented all 
algorithms in the computer algebra system \textsc{Singular}, and compare the different approaches with respect to timings and memory usage. The test examples 
originate from questions in algebraic geometry, where the use of the Jacobian 
criterion is impractical due to the number and size of the minors involved.
\end{abstract}

\section{Introduction}

In classical algebraic geometry, explicit constructions of algebraic varieties 
with prescribed properties play an important role, for example, for existence 
and unirationality results in moduli problems, see e.g. \cite{NP1,Sch}. In 
many situations, the aim is to construct a {\em smooth} variety satisfying 
certain properties. For example, when considering a family of algebraic curves with a given
Hilbert polynomial, the arithmetic genus of the curve is determined. 
The geometric genus, however, differs from the arithmetic 
genus by the delta invariant, which measures the singularities of the curve. 
Hence the presence of singularities affects the geometric genus, leading to a 
different topological type of the curve. Passing to the next dimension, recall that
algebraic surfaces have been classified in 
the Enriques-Kodaira classification, see e.g. \cite{BHPV}. Especially for 
surfaces of general type, a full understanding of the moduli spaces and 
explicit constructions of canonical rings are still lacking. Although there 
are various techniques to construct surfaces with prescribed invariants, see for example
\cite{LP,NP1}, the constraint of smoothness often requires testing and, in practice, 
this turns out to be a fundamental obstacle. It is thus of equal importance 
to have a test which is both fast in determining smoothness and non-smoothness.

The standard method for testing smoothness is the Jacobian criterion, (cf. any 
textbook on computational algebraic geometry, e.g. \cite[Cor. 5.6.14]{GP}). 
Given an affine algebraic variety $X=V(I)\subseteq \mathbb{A}^n$ of 
codimension $c$ defined by an ideal $I=\langle f_1,\ldots,f_s \rangle\subseteq 
k[x_1,\ldots,x_n]$ over an algebraically closed field $k$, we compute the 
dimension of the vanishing locus of the Jacobian ideal $J$, which is generated by
$c\times c-$minors of the Jacobian matrix $(\frac{\partial f_i}{\partial x_j})$
on $X$. This can be done by computing a Gr\"obner basis of the ideal $I+J$. 
However, the number $\binom{n}{c}\cdot\binom{s}{c}$ of minors can be very 
large, and the Gr\"obner basis determines the complete scheme structure of the 
singular locus of $X$, which is not required to check smoothness. As a 
result, this approach will be rather inefficient and often even impractical.

In this paper, we describe an algorithm for determining smoothness, which is 
based on an idea from Hironaka's famous proof of resolution of singularities. 
The (implicitly stated) termination criterion provides a smoothness criterion, 
which does not require computation of the $c\times c-$minors of the Jacobian 
matrix. The key idea behind this smoothness test is the fact that each 
non-singular variety is locally a complete intersection. That is, it can be
covered by suitable open subsets, on each of which we truly see a complete 
intersection. Such a covering can be computed without too much effort and, in 
relevant examples, our approach significantly extends the limits of 
practicability of the Jacobian criterion.  In addition, the construction is 
inherently parallel with only minimal communication overhead. In fact, 
application of the criterion 
can in some cases be faster than computing a single minor of the 
Jacobian matrix. Our algorithm is implemented in the library 
\texttt{smoothtst.lib} \cite{smoothtst} for the 
computer algebra system \textsc{Singular}.

This paper is organized as follows: In Section \ref{sec Hironaka}, we extract 
the relevant part of Hironaka's smoothness criterion from the constructive, 
simplified desingularization approach of \cite{BEV}. In Section 
\ref{sec algorithm1}, we turn this criterion into an algorithm and then, in Section
\ref{sec hybrid}, refine 
it into a hybrid approach, which combines it with the use of the Jacobian criterion
in smaller codimension. 
Section \ref{sec examples} introduces some settings and constructions from 
algebraic geometry where an efficient smoothness test is required, and presents 
explicit examples thereof. These are then used in Section \ref{sec timings} 
to compare our new approaches with the standard technique based on the Jacobian 
criterion.

\section{Hironaka's Smoothness Criterion}\label{sec Hironaka}

In 1964, Hironaka proved existence of resolutions of singularities in 
characteristic zero  \cite{Hir}. 
He introduced standard bases to achieve this goal. Though his 
proof is non-constructive in certain parts, standard bases are by no means the
only algorithmic considerations introduced there. The termination criterion, 
which he uses, provides a smoothness criterion that does not involve the 
computation of the ideal of codimension-sized minors of the Jacobian matrix.\\

For this article, we assume $k$ to be an algebraically closed field of 
characteristic zero unless explicitly stated otherwise. The general line of 
arguments is still valid, if we drop the condition on the field to be 
algebraically closed, but everything needs to be stated with significantly 
more care.

\begin{df}\cite{Hir} 
Let $(X,p) \subset ({\mathbb A}_k^n,p)$ be a germ with defining ideal 
$I_{X,p} \subset k\{\underline{x}\}:=k\{x_1,\ldots,x_n\}$ generated by 
 $f_1,\dots,f_s$, and assume that these power series form a standard basis 
of $I_{X,p}$ with respect to some local degree ordering. Assume further 
that the power series $f_i$ are sorted by increasing order. We then denote 
by $\nu^*(X,p)$ the sequence of orders of the $f_i$.
\end{df}

Making use of the sequence $\nu^*(X,p)$, we can formulate the termination 
condition in Hironaka's proof in the form of the following lemma, which is 
already implicitly present and used in \cite{Hir}:

\begin{lem} The germ
$(X,p) \subset ({\mathbb A}_k^n,p)$ is singular at $p$ if and only if
$$\nu^*(X,p) >_{lex} (\underbrace{1,\dots,1}_{\operatorname{codim}(X)})$$ 
with respect to the lexicographical ordering $>_{lex}$.
\end{lem}

Of course, it is impractical to consider each point of a given variety 
separately, however this is a general problem when translating Hironaka's 
ideas into an algorithmic approach. Several groups have successfully tackled 
the problem of making Hironaka's approach completely constructive since the 
end of the 1980s (cf. \cite{BM}, \cite{BEV}, \cite{EH}), and have
created algorithmic approaches to desingularization. We will roughly follow here
the approach of Bravo, Encinas and Villamayor, drop the condition of 
being a standard basis, and use the notion of the order of an ideal at a 
specific point $p$, which is precisely the first entry of $\nu^*(X,p)$. 
To make this notion accessible to computations, we state everything and do 
all calculations in affine charts, that is from the algebraic point of 
view, in polynomial rings. For the more general situation of an algebraic 
variety, the results can then be applied by means of a suitable cover of the 
variety by affine charts. 

\begin{df}
Let $W \subseteq {\mathbb A}_k^n$ be a smooth variety and let 
$I=\langle f_1,\dots,f_s \rangle \subseteq k[\underline{x}] $ be a 
radical ideal defining a variety $X \subseteq W$. Then
$$\operatorname{ord}_p(I):= 
       \max \{m \in {\mathbb N} \mid I \subseteq {\mathfrak m}_{W,p}^m\}.$$
\end{df}

\begin{rem}
Given an ideal $I$ as above and adapting \cite{BEV} to our setting, we define 
the ideal sheaf $\Delta(I)$ as the ideal sheaf locally generated by 
$f_1,\dots,f_s$ and all partial derivatives thereof with respect to a regular 
system of parameters of $W$. 
Then the locus of order at least two is precisely the vanishing locus of 
$\Delta(I)$, that is, the order of the ideal is (at most) one everywhere, if 
and only if $1 \in \Delta(I)_w$ for all $w \in W$.  
\end{rem}

Beware of a misunderstanding here: The above remark does not mean
 that there is one of the generators of $I$ which defines a smooth 
hypersurface, nor does it mean that we can use the same regular system 
of parameters everywhere. It only implies that, at each point of $X$, 
there is some generator defining a smooth hypersurface $W_1 \subset W$ in 
some open neighborhood of the point (for an algorithmic treatment of the 
general case see \cite[Section 4.2]{FK1}).

\begin{rem}
Covering $X$ by affine charts corresponding to the complements of 
hypersurfaces, each containing the singular locus of one of the $f_i$, we can 
locally use $W \cap V(f_i)$ as the new ambient space $W_1$ and iterate our 
considerations\footnote{For those readers who are familiar with algorithmic
resolution of singularities, we would like to add two remarks: Due to 
the low order of the ideal and the absence of exceptional divisors from 
previous blow-ups, this is a specifically simple instance of 
the local existence of a hypersurface of maximal contact and a descent in 
ambient dimension by means of a coefficient ideal, see \cite{BEV} or 
\cite{FK1} for an implementable general construction.}. If we ever encounter 
a maximal order exceeding one during this iteration -- in any of the charts -- 
we know that $X$ is singular.

The key idea behind the smoothness test is the fact that each non-singular 
variety is locally a complete intersection. Passing to a covering such that 
the variety is a complete intersection in each chart, all generators 
of the corresponding ideal are of order one.
\end{rem}

The real difficulty in the 
process is the task of determining a suitable covering. Actually, the crucial 
point in turning this strategy into an algorithm is the choice of suitable 
regular systems of parameters; luckily this has already been made algorithmic 
in the more general case of the resolution of singularities in \cite{FK1}, \cite{FK4}, and \cite{reslib}. Given such a system of parameters, it is then straight forward 
to perform the iteration and, hence, to test for smoothness. We give the 
details in the subsequent section.

\section{An Algorithmic Smoothness Test}\label{sec algorithm1}
For the smoothness test described in this section, we assume that the input 
is an affine variety $X \subseteq W\subseteq {\mathbb A}_k^n$ 
such that the ambient variety $W$ is ${\mathbb A}_k^n$ at the start of the algorithm, 
that is, $I_W=\langle 0 \rangle \subset k[x_1,\dots,x_n]$. If we 
are dealing with more general varieties, the smoothness test needs to be 
applied separately to each affine chart of the covering of the variety. In particular, 
a subvariety of projective space can be covered by the standard affine charts.

The top-level algorithm for our smoothness test is a rather simple recursion, for 
the details see Algorithm \ref{alg smoothnesstest}. 
This algorithm is just a reformulation of the ideas of Section \ref{sec Hironaka} and, 
of course, nothing but a hollow shell without the two workhorses described 
in Algorithms \ref{alg deltacheck} and \ref{alg descendembsmooth}, which we 
will now discuss in detail.

\begin{algorithm}[h]
\caption{Smoothness Test for Subvarieties}
\label{alg smoothnesstest}
\begin{algorithmic}[1]

\REQUIRE       Ideals    $I_W= \langle g_1,\dots,g_r \rangle,\, I_X= \langle f_1,\dots,f_s \rangle \subseteq 
                           k[x_1,\dots,x_n]$ of affine varieties $W$, $X$, respectively,
            and
             $g\in k[x_1,\dots,x_n]$ such that 
             
             \begin{itemize}[leftmargin=*]
             \item $I_W \subseteq I_X$, and
             \item $W$ is a non-singular complete intersection on $D(g)$.
             \end{itemize}
\ENSURE {\texttt true} if $X$ is non-singular, {\texttt false} otherwise.

\IF{ $I_W =I_X$ on $D(g)$}
    \RETURN {\texttt true}
\ENDIF    
\IF {{\bf not} $\operatorname{DeltaCheck}(I_W,I_X,g)$} 
   \RETURN {\texttt false}
\ENDIF\label{Jac1}
\STATE $L:=\operatorname{DescendEmbeddingSmooth}(I_W,I_X,g)$\label{Jac2}
\FORALL {($I_U,I_{X\mid U},g_U)\in L$}
     \IF {{\bf not} $\operatorname{SmoothnessTest}(I_U,I_{X|U},g_U)$ }
         \RETURN {\texttt false}
    \ENDIF
\ENDFOR    
\RETURN {\texttt true}

\end{algorithmic}
\end{algorithm}

\begin{algorithm}[h]
\caption{DeltaCheck}
\label{alg deltacheck}
\begin{algorithmic}[1]

\REQUIRE         Ideals    $I_W= \langle g_1,\dots,g_r \rangle,\, I_X= \langle f_1,\dots,f_s \rangle \subseteq 
                           k[x_1,\dots,x_n]$ of affine varieties $W$, $X$, respectively,
            and
             $g\in k[x_1,\dots,x_n]$ such that 
             
             \begin{itemize}[leftmargin=*]
             \item $I_W \subseteq I_X$, 
             \item $W$ is a non-singular complete intersection on $D(g)$.
             \end{itemize}

\ENSURE  {\texttt true} if $1 \in \Delta(I_X)$ on $D(g)$, {\texttt false} otherwise.

\smallskip\noindent\emph{Boundary case: 
$W={\mathbb C}^n$, $g=1$, that is,
               $x_1,\dots,x_n$ is regular system:}\smallskip

\IF {$I_W=\langle 0\rangle$ {\bf and }$g=1$}
    \IF {$1 \in \langle f_1,\dots,f_s, \frac{\partial f_1}{\partial x_1},\dots, \frac{\partial f_s}{\partial x_k}\rangle$}
      \RETURN {\texttt true}
    \ELSE 
       \RETURN {\texttt false}
     \ENDIF
\ENDIF

\smallskip\noindent\emph{Initialization:}\smallskip

  \STATE $Q:=\langle 0\rangle$
   \STATE $L_1:=\{ r\times r-\text{submatrices } M\text{ of }\operatorname{Jac}(I_W) \text{ with }\operatorname{det}(M)\neq0\}$

\smallskip\noindent\emph{Main Loop: Covering by complements of the minors:}\smallskip

    \WHILE {$L_1\neq \emptyset$ {\bf and} $g \not\in Q$}
      \STATE choose $M \in L_1$
      \STATE $L_1:= L_1 \setminus \{M\}$
      \STATE $q:=\operatorname{det}(M)$
      \STATE $Q:=Q+\langle q\rangle$
      \STATE compute an $r\times r-$matrix $A$ such that
      $$A \cdot M = q \cdot \operatorname{Id}_{r}$$ 

\smallskip\noindent\emph{Determine the components of $\Delta(J)$ not lying inside $V(q) \cup V(g)$:}\smallskip

            \STATE $C_M:= I_X \; +  \;
             \left\{q \cdot \frac{\partial f_i}{\partial x_j}
             -\sum\limits_{k \text{ column of }M \atop 
                    l\text{ row of } M} 
              \frac{\partial g_l}{\partial x_j} A_{lk} 
              \frac{\partial f_i}{\partial x_k} \left| 
              \begin{footnotesize}
                \begin{array}{c}
                 1 \leq i \leq s,\\
                 j \text{ not a column of } M
                \end{array}
              \end{footnotesize}
              \right. \right\}$\label{line derivatives}
     \IF {$q\cdot g \not\in \sqrt{C_M}$}
         \RETURN \texttt false
     \ENDIF
   \ENDWHILE

\RETURN {\texttt true}
\end{algorithmic}
\end{algorithm}

\begin{notation}
For stating the algorithms of this section, we adhere to standard notation of 
algebraic geometry denoting the complement of the vanishing locus of a 
polynomial $g \in k[x_1,\dots,x_n]$  by $D(g)$. The equality test for 
$I_W \subseteq I_X$ on $D(g)$ then amounts to checking whether $g \cdot I_X$ is
contained in $\sqrt{I_W}=I_W$ (where the last equality holds, because $W$ is 
smooth). 
\end{notation}

Algorithm \ref{alg deltacheck}, {\it DeltaCheck}, is designed to detect non-singularity as soon 
as possible with a minimal amount of time and memory consumption. To this end, 
it first treats the 'lucky' trivial case, where the ambient space $W$ is just 
${\mathbb A}_k^n = \operatorname{Spec}(k[x_1,\dots,x_n])$, and the
regular system of parameters is hence given by the variables. If a regular
system of parameters needs to be determined, we can make use of the 
smoothness of the ambient space $W$, which implies that at each point of 
$W \cap D(g)$ at least one of the $r \times r$ minors of the Jacobian 
matrix of $I_W=\langle g_1,\dots,g_r \rangle$ is non-zero. Thus, the 
complements of these minors give rise to an open covering of $W \cap D(g)$, and,
on each open set, the columns which are unused with respect to the invertible minor 
lead to a regular system of parameters. But, since coordinate changes are 
expensive, the explicit passage to the new coordinates (and back) is avoided 
by an appropriate modification of the derivatives (see line \ref{line derivatives}). Eventually, we only
need to check whether $qg \in \sqrt{C_M}$, because this amounts to the same
as checking whether $V(C_M)$ is contained in the union of the
hypersurfaces $V(g)$ and $V(q)$. 

\begin{rem}
It is, in general, not a good idea to run this computation for all 
$\binom{n}{r}$ minors. Instead, the computation should stop as soon as all of
$W \cap D(g)$ is covered by the open sets. This explains the second 
termination condition of the main loop. As an enhancement, it is 
possible to further reduce the number of charts by expressing 
$g$ in terms of the minors of $\operatorname{Jac}(I_W)$ and then 
to only consider those minors appearing.\\
Expressing the element $g$ in terms of the generators of an ideal can easily 
be implemented in {\sc Singular}: the command {\tt lift} performs this 
operation on the basis of a single standard basis computation. Thus the 
previously considered enhancement is actually a trade-off, as it avoids 
redundant calculations at the cost of an additional standard basis computation.
\end{rem}

\begin{algorithm}[h]
\caption{DescendEmbeddingSmooth}
\label{alg descendembsmooth}
\begin{algorithmic}[1]
\REQUIRE $(I_W,I_X,g)$ where $I_X=\langle f_1,\dots,f_s \rangle$, 
         $I_W=\langle g_1,\dots,g_r \rangle \subseteq k[x_1,\dots,x_n]$
             are ideals  of affine varieties $W$, $X$, respectively, and
             $g\in k[x_1,\dots,x_n]$ such that
             \begin{itemize}[leftmargin=*]
             \item $I_W \subseteq I_X$,
             \item $W$ a non-singular complete intersection on $D(g)$,  and
             \item the order of $I_X$ is at most $1$ at each point of $W$
             \end{itemize}
\ENSURE set of triples $(I_{Z_i},I_{X|Z_i},D(g_i))$ such that
              \begin{itemize}[leftmargin=*]
              \item $Z_i$ contains the support of $I_X$ in $D(g_i)$, and
                    is a non-singular complete intersection on $D(g_i)$
              \item $I_{X|Z_i}$ is the restriction
                    of $I_X$ to $Z_i \cap D(g_i)$
              \end{itemize}

\smallskip\noindent\emph{Check whether we can avoid an open covering:}
\smallskip

     \IF {$V(f_i)$ non-singular hypersurface in $W \cap D(g)$ for some $i$}\label{line choose fi}
        \STATE $I_Z:=\langle g_1,\dots,g_r,f_i\rangle$
        \STATE $I_{X|Z}:=I_X$
        \RETURN \{($I_Z$,$I_{X|Z}$,$g$)\}\label{line end choose fi}
     \ENDIF
\smallskip     \noindent\emph{Find an open covering using that $\bigcap_{i=1}^s Sing(f_i) =\emptyset$:}\smallskip
      
      \STATE Express $g$ in terms of ideals of singular loci and
             find an open covering of $W \cap D(g)$ by complements of 
             hypersurfaces such that
      \begin{itemize}
      \item $I_{Z_j}=\langle g_1,\dots,g_r,f_{i_j} \rangle$
      \item $I_{X|Z_j}=I_X$
      \item $g_j=g\cdot h_j$ where $\operatorname{Sing}(f_{i_j}) \subseteq V(h_j)$ and
            $1 \in \langle h_1,\dots,h_k \rangle$
      \end{itemize}
\STATE $L:=\{(I_{Z_1},I_{X|Z_1},g_1),\dots,(I_{Z_k},I_{X|Z_k},g_k)\}$
\RETURN $L$

\end{algorithmic}
\end{algorithm}

The general idea of the Algorithm \ref{alg descendembsmooth}
, 
{\it DescendEmbeddingSmooth}, is the one described in
its lines \ref{line choose fi} to \ref{line end choose fi}: Find an element of $I_X$, which defines a non-singular 
hypersurface in $W \cap D(g)$, use this hypersurface as a new ambient space 
and return the new triple. Unfortunately, the situation is more complicated, 
as, in general, such a hypersurface does not exist globally.
An open covering is thus required to find such a hypersurface on each open 
subset (as already outlined in the previous section). On the other hand, we 
know that the maximal order of $I_X$ at each point of $W \cap D(g)$ is at 
most one, which implies that the intersection of the singular loci of the 
hypersurfaces $V(f_i) \cap D(g)$ in $W \cap D(g)$ is empty. When expressing $g$,
which of course is a unit on $D(g)$, in terms of the generators of the ideals
of the singular loci, only some of the generators actually appear in the 
corresponding $k[x_1,\dots,x_n]$-linear combination with non-zero coefficients. 
These are the polynomials, for which the complements of the vanishing loci
define the desired open covering. As soon as this covering has been found,
we can proceed as in steps \ref{line choose fi} to \ref{line end choose fi} for each of the open sets.\\

Note that, to apply Algorithm \ref{alg descendembsmooth} in Algorithm 
\ref{alg smoothnesstest}, we actually need that
$Z_i$ contains the locus of maximal order of $I_X$ on $D(g)$. But at this point
we already know that the maximal order is one, whence the locus of maximal 
order is just the vanishing locus of $I_X$.
Note also that $I_{X|Z_i}$ is the
coefficient ideal of $I_X$ with respect to $Z_i$ since the maximal order of $I_X$ is one.

If we assume that Algorithm \ref{alg smoothnesstest} was initially called with
$W= {\mathbb A}_k^n$, the construction  at each level of descent in ambient dimension provides an ambient space which is a 
complete intersection on $D(g)$.
If the top-level ambient space, that is, the original input, is not a complete 
intersection, it may be necessary to find an open covering such that this 
smooth ambient space, which locally is a complete intersection, is specified 
by a complete intersection on each open set itself, before the smoothness
test for $X \subseteq W$ can be applied. For an affine variety $W$,
this can be achieved by applying Algorithm \ref{alg descendembsmooth} to $W\subseteq {\mathbb A}_k^n$ codimension times.

\begin{rem}
The above Algorithms \ref{alg smoothnesstest} and \ref{alg deltacheck} are 
inherently parallel, as both rely on open coverings and the computations 
in the respective open sets are independent of each other. More precisely,
the main loop of Algorithm \ref{alg smoothnesstest}, in which the recursion occurs, can
be distributed to multiple cores and the computation of $\Delta(I_X)$ on the 
different sets $D(q) \cap D(g)$ also readily allows parallelization.

Due to the boolean return type, the total communication overhead of 
parallelizing one of these steps originates from handing data to the routine, 
not from returning a result. The embedding dimension, on the other hand, 
provides an upper bound for the depth of the recursion. In combination, this 
implies that the communication overhead for a fully parallel approach will be
insignificant compared to the benefits in computation time. 
\end{rem}

\begin{rem}
For practical purposes, in line \ref{line choose fi} of Algorithm 
\ref{alg descendembsmooth}, it is worth to try a general linear combination 
of $f_1,\dots,f_s$ as well, before passing to an open cover.
\end{rem}

\begin{rem}
There is one potential drawback to the approach outlined in this section: In 
comparison to the standard approach of using the Jacobian criterion, the 
presented algorithms do not only rely on arithmetic operations on a fixed set 
of polynomials and one subsequent standard basis computation. They need, for example,
standard bases computations for the radical membership tests. This does not 
only make them more sensitive to intermediate coefficient swell, but also 
introduces imponderabilities with regard to the complexity of a given example due 
to the well known difference between average and worst-case complexity of computing 
standard bases.
\end{rem}

\begin{rem}
Up to this point, we always assumed that $\operatorname{char} k=0$. This may 
seem to be an obvious restriction, given that we are extracting the approach from the proof 
of resolution of singularities, which is a famous open problem in 
positive characteristic. 
Considering the obstacles for generalizing Hironaka's approach to positive 
characteristic a bit more closely (as was done, e.g., in \cite{Ha1}, and, with 
focus on constructing examples, also in \cite{FK3}), it is 
obvious that all of the unresolved complications arise during the sequence 
of blow-ups, not within the smoothness-test. So this issue does not pose a serious
obstruction to using the above algorithms over a perfect field in positive 
characteristic in order to check smoothness of a reduced variety $X \subseteq {\mathbb A}_k^n$. 
One could, thus, also envision the use of a modular approach to the characteristic zero problem.
\end{rem}

\section{A Hybrid Approach\label{sec hybrid}}

The main drawback of the classical smoothness test by means of the Jacobian 
criterion is the high number of minors to be computed. It may therefore not 
only turn out to be time consuming, but also requires huge amounts of memory. 
In practice, the latter is the more restricting aspect, as it may cause a 
premature termination of the computation due to lack of memory either while
computing the minors or while determining a Gr\"obner basis of the Jacobian ideal. The new 
approach can also lead to involved Gr\"obner basis 
computations -- in particular, if the polynomials defining the open sets $D(g)$
are of high degree and contain many terms. 

So both approaches have rather obvious drawbacks, but the one of the classical
approach is most relevant if the codimension of the given variety is large and
the variety is far from being a complete intersection. For the new approach, 
however, the main drawback does not become significant, before several 
descents in ambient dimension have been applied. Each descent in ambient 
dimension lowers the codimension of the variety. Hence, it makes sense to use
the new algorithm as a preparation step for applying the Jacobian criterion in a smaller codimension.

In Algorithm \ref{alg smoothnesstest}, this hybrid approach can easily be 
introduced by adding a new input parameter controlling the depth of the descent
in ambient dimension and checking between lines \ref{Jac1} and \ref{Jac2}, whether the maximal
desired number of descents has been made. If so, we call Algorithm
\ref{alg embJac} to check smoothness by the Jacobian criterion and return its
result.

\begin{algorithm}[h]
\caption{Jacobian criterion for $X \subseteq W$}
\label{alg embJac}
\begin{algorithmic}[1]

\REQUIRE   Ideals    $I_W= \langle g_1,\dots,g_r \rangle,\, 
                     I_X= \langle f_1,\dots,f_s \rangle \subseteq 
                           k[x_1,\dots,x_n]$ of affine varieties $W$, $X$, 
                     respectively, and
             $g\in k[x_1,\dots,x_n]$ such that 
             
             \begin{itemize}[leftmargin=*]
             \item $I_W \subseteq I_X$, and
             \item $W$ is a non-singular complete intersection on $D(g)$.
             \end{itemize}
\ENSURE {\texttt true} if $X$ is non-singular, {\texttt false} otherwise.

\smallskip

  \STATE $Q:=\langle 0\rangle$
   \STATE $L_1:=\{ r\times r-\text{submatrices } M\text{ of }\operatorname{Jac}(I_W) \text{ with }\operatorname{det}(M)\neq0\}$

\smallskip\noindent\emph{Covering by complements of the minors:}\smallskip

    \WHILE {$L_1\neq \emptyset$ {\bf and} $g \not\in Q$}
      \STATE choose $M \in L_1$
      \STATE $L_1:= L_1 \setminus \{M\}$
      \STATE $q:=\operatorname{det}(M)$
      \STATE $Q:=Q+\langle q\rangle$
      \STATE compute an $r\times r$-matrix $A$ such that
      $$A \cdot M = q \cdot \operatorname{Id}_{r}$$ 

\smallskip\emph{Determine Jacobian matrix of $I_X$ w.r.t. system of parameters for $W$:}\smallskip

            \STATE $Jac :=
             {\left(q \cdot \frac{\partial f_i}{\partial x_j}
             -\sum\limits_{k \text{ a column of }M \atop 
                    l\text{ a row of } M} 
              \frac{\partial g_l}{\partial x_j} A_{l,k} 
              \frac{\partial f_i}{\partial x_k} \right)_{{1\leq i\leq n} \atop
                       {j \text{ not a column of } M}}}$ \medskip
     \STATE $J := I_X + \left\langle (\dim(W)-\dim(X)) \text{-sized minors of }Jac \right\rangle$
     \IF {$q\cdot g \not\in \sqrt{J}$}
             \RETURN \texttt false
     \ENDIF
   \ENDWHILE

\RETURN {\texttt true}
\end{algorithmic}
\end{algorithm}

Algorithm \ref{alg embJac}, which applies the Jacobian criterion in
a relative situation, is again a hybrid between two algorithms that we have seen 
before. The computation of a regular system of parameters (possibly after
passing to an open covering) and the computation of the Jacobian matrix 
with respect to this system coincide
precisely with the respective steps of Algorithm \ref{alg deltacheck}.
The remaining steps are a straight forward application of the Jacobian criterion to the
before-computed Jacobian matrix.

This hybrid approach inherits the inherently parallel structure of the new
algorithm, lowers the combinatorial complexity of the Jacobian criterion by
applying several descents of ambient dimension beforehand, and avoids -- to a 
certain extent -- the imponderabilities of the large number of radical 
membership tests arising from covering and descending as many times as the 
codimension requires.
 In this way, we obtain a viable 
divide-and-conquer approach for parallelizing the Jacobian criterion. 
      
\section{Applications in Algebraic Geometry\label{sec examples}}

In this section, we describe several constructions from algebraic geometry, which
all aim at obtaining varieties with specific properties. All resulting ideals 
are specified as examples in our library \texttt{smoothtst.lib}. We use these
examples in Section \ref{sec timings} to study the performance of our
smoothness test in realistic usage scenarios. 

The constructions use the technique of Kustin-Miller unprojection to produce
sequences of birationally equivalent varieties. Unprojection theory was
introduced by A. Kustin, M. Miller, M. Reid, and S. Papadakis (see \cite{KM,
R1, PR, P1}) to construct Gorenstein rings of high codimension and to act as a
substitute for a structure theorem for Gorenstein rings of codimension $\geq4$. We give a short
outline of the fundamental construction used in the examples, for details see
\cite{BP3}. This construction is implemented in the package
\texttt{KustinMiller} \cite{BPlib} for the computer algebra system
\textsc{Macaulay2} \cite{GS}.

Suppose that $R$ is a positively graded polynomial ring over a field. Given
two homogeneous ideals $I\subset J\subset R$ defining Gorenstein rings
$R/I$ and $R/J$ such that $\dim R/I=\dim R/J+1$, that is, $V(J)$ is
codimension $1$ subvariety of $V(I)$, we will construct a new Gorenstein ring,
the unprojection ring. Geometrically the construction corresponds to the
contraction of the subvariety $V(J)\subset V(I)$. By \cite[Proposition~3.6.11]%
{BH} there are $k_{1},k_{2}\in\mathbb{Z}$ such that $\omega_{R/I}=R/I(k_{1})$
and $\omega_{R/J}=R/J(k_{2})$. Suppose that $k_{1}>k_{2}$. This implies that
the unprojection ring is also positively graded.

\begin{df}
\label{def unproj ring}\cite{PR} Let $\iota:J/I\rightarrow R/I$ be the
inclusion morphism and let $\phi:J/I\rightarrow R/I$ be a homomorphism of degree
$k_{1}-k_{2}$ such that $\operatorname{Hom}_{R/I}(J/I,R/I)$ is generated as an
$R/I$-module by $\iota$ and $\phi$. The Kustin--Miller unprojection
ring defined by $\phi$ is the graded algebra
$$(R/I)[T]/U\quad\text{ where }\quad U = \left\langle Tu-\phi(u)\mid u\in J/I\right\rangle,$$ and $\deg(T):=k_{1}-k_{2}$.
\end{df}

Note that $(R/I)[T]/U\cong R[T]/\widetilde U$, where $\widetilde U$ is the inverse image of $U$
under the natural map $R[T]\rightarrow(R/I)[T]$.

\begin{prop}
\cite{KM,PR} The $R$-algebra $R[T]/U$ is Gorenstein and, up to isomorphism,
independent of the choice of $\phi$.
\end{prop}

\begin{example}
\label{ex del pezzo} Consider a cubic surface in $\mathbb{P}^3$ with six 
disjoint lines. Such a surface corresponds to the blow-up of $\mathbb{P}^2$ 
in six points (for an implementation of this construction, see the \textsc{Macaulay2} package 
\texttt{clebsch.m2} \cite{cubic}). Applying the Kustin-Miller unprojection 
construction to the ideal $J$ of such a line corresponds to the Castelnuovo 
blow-down of the line and results in a Del Pezzo surface in $\mathbb{P}^4$ of 
degree $4$. By iterating this construction, we obtain a sequence of Del Pezzo 
surfaces of degree $2+c$ and codimension $c$ in $\mathbb{P}^{2+c}$ for 
$c=1,\ldots,7$. See \cite[Section 4.1]{NP2} for more details. \\
\end{example}

\begin{example}
\label{ex veronese} Consider the codimension $2$ complete intersection threefold $X_2\subset \mathbb{P}^5$ defined by 
$$I=\langle  x_1x_3-y_1y_2,\ x_2x_4- y_1y_2\rangle$$ and the ideals 
$J_1=\langle y_1,x_1,x_2 \rangle$, 
$J_2=\langle y_2,x_2,x_3 \rangle$, 
$J_3=\langle y_1,x_3,x_4 \rangle$, and 
$J_4=\langle y_1,x_4,x_1 \rangle$. 
Applying the Kustin-Miller construction iteratively for the ideal $J_1$ 
introducing the new variable $T_1$, for $J_2+\langle T_1 \rangle$ introducing 
$T_2$, for $J_3+\langle T_1,T_2 \rangle$ introducing $T_3$, and for 
$J_4+\langle T_1,T_2,T_3 \rangle$ introducing $T_4$, we obtain a sequence 
$X_2,\ldots,X_6$ of threefolds of codimensions $2,\ldots,6$. 
Note that $X_2,X_3,X_4$ are singular, whereas $X_5$ and $X_6$ are smooth.
\end{example}

\begin{example}\label{ex cyclic}
Let $R=k[  x_{1},...,x_{n}]  $ be a polynomial ring over a field $k$ and let
$I_{d}\left(  R\right)  $ be the Stanley-Reisner ideal of the boundary complex of
the cyclic polytope of dimension $d$ with vertices $x_{1},\ldots,x_{n}$. (Recall 
that the Stanley-Reisner ideal is the monomial ideal generated by the 
non-faces of the complex.)
As shown in \cite{BP2}, Kustin-Miller unprojection yields a recursion
for the ideals $I_{d}\left(  R\right)  $: For $d$ even, apply the construction
with $T=x_{n}$ to $I=I_{d}\left(  k\left[  x_{1},...,x_{n-1}\right]  \right)  $
and $J=I_{d-2}\left(  k\left[  z,x_{2},...,x_{n-2}\right]  \right)  $ 
considered as ideals in $k\left[  z,x_{1},...,x_{n-1}\right]  $ and quotient 
by $\langle z\rangle  $. For $d$ odd one can proceed in a similar way, for details see \cite[Section 4]{BP2}.
Note that the varieties defined by the ideals $I_d(R)$ are not smooth.
\end{example}

\section{Efficiency and Timings}

\label{sec timings}

The Hironaka style smoothness test as presented in Section \ref{sec algorithm1} together with the hybrid variant discussed in Section \ref{sec hybrid} are implemented in the
\textsc{{Singular}} library {\texttt{smoothtst.lib}}, see \cite{smoothtst}. In
this section, we compare the performance of these approaches to the smoothness
test given by the Jacobian criterion (see Table \ref{tab timings}, columns \emph{smoothtst}, \emph{hybrid}, and \emph{Jacobian}). 

Note that, when computing minors of 
the Jacobian matrix of $X=V(I)\subset\mathbb{A}^n$, we consider all results 
of arithmetic operations in $k[x_1,\dots,x_n]/I$, that is, we reduce modulo $I$. 
This leads to a significant improvement of the performance of the Jacobian 
criterion. 
In the hybrid approach, we exemplarily give timings for descending to codimension $2$,
that is, applying the Jacobian criterion for $2\times 2-$minors.


As discussed above, the communication overhead of a parallel version of the 
algorithm depends only on the depth of the recursion tree generated in the 
course of the algorithm. This depth is bounded by the embedding dimension. 
Moreover, the return type of the recursion is of type boolean. As a result, 
the algorithm has basically zero communication overhead. Hence, to show
the full potential of the algorithm with respect to massively parallel computations, we give
simulated parallel timings even if the number of processes exceeds the number
of cores available on our machine (see subcolumns \emph{time} for sequential and \emph{parallel} for parallel timings). 

All timings are in seconds on an Intel Xeon E5-2690 machine with $32$ cores, $2.9$ GHz, and $192$ GB of
RAM running a Linux operating system. 
In addition to timings, we also do a comparison with respect to memory 
usage by indicating the maximum amount of memory used by \textsc{Singular} in 
megabytes when running the algorithms on a single core (subcolumns  \emph{mem}).
Dashes
indicates that the respective computation 
did not finish within $10000$ seconds or 
used more than $20$GB of RAM.

In the table we also
indicate, whether the specific example is smooth (column \emph{smooth}).  We have verified that for
all examples, where the Jacobian criterion finishes or where the 
presence/absence of singularities is a known from theoretical considerations, our implementation of the new algorithm comes 
to the same result. 

\begin{table}[h]
\caption{Timings and Memory Usage}%
\begin{tabular}
[c]{l|c|lll|ccc|cl|}
&  &  &  & \hspace{-18mm}smoothtst &  &  & \hspace{-30mm}hybrid &  &
\hspace{-6mm}Jacobian\\
& smooth & time & \multicolumn{1}{|l}{parallel} & \multicolumn{1}{|l|}{mem} &
time & \multicolumn{1}{|c}{parallel} & \multicolumn{1}{|c|}{mem} &
\multicolumn{1}{|l|}{time} & \multicolumn{1}{|c|}{mem}\\\hline\hline
$\mathcal{I}_{1}(6)$ & yes & \multicolumn{1}{|c}{0.24} &
\multicolumn{1}{|c}{0.07} & \multicolumn{1}{|c|}{0.22} & 0.18 &
\multicolumn{1}{|c}{0.05} & \multicolumn{1}{|c|}{0.22} & 2.5 &
\multicolumn{1}{|c|}{34}\\\hline
$\mathcal{I}_{1}(7)$ & yes & \multicolumn{1}{|c}{0.60} &
\multicolumn{1}{|c}{0.17} & \multicolumn{1}{|c|}{0.24} & 0.35 &
\multicolumn{1}{|c}{0.10} & \multicolumn{1}{|c|}{0.22} & 310 &
\multicolumn{1}{|c|}{1300}\\\hline
$\mathcal{I}_{1}(8)$ & yes & \multicolumn{1}{|c}{0.86} &
\multicolumn{1}{|c}{0.22} & \multicolumn{1}{|c|}{0.32} & 0.64 &
\multicolumn{1}{|c}{0.15} & \multicolumn{1}{|c|}{0.23} &
\multicolumn{1}{|c|}{$-$} & \multicolumn{1}{|c|}{$>20000$}\\\hline\hline
$\mathcal{I}_{2}(3)$ & yes & \multicolumn{1}{|c}{0.22} &
\multicolumn{1}{|c}{0.04} & \multicolumn{1}{|c|}{0.14} & 0.08 &
\multicolumn{1}{|c}{0.02} & \multicolumn{1}{|c|}{0.14} & 0.05 &
\multicolumn{1}{|c|}{4.2}\\\hline
$\mathcal{I}_{2}(4)$ & yes & \multicolumn{1}{|c}{160} &
\multicolumn{1}{|c}{9.1} & \multicolumn{1}{|c|}{27} & 40 &
\multicolumn{1}{|c}{4.9} & \multicolumn{1}{|c|}{190} & 15 &
\multicolumn{1}{|c|}{450}\\\hline
$\mathcal{I}_{2}(5)$ & yes & \multicolumn{1}{|c}{$-$} &
\multicolumn{1}{|c}{$-$} & \multicolumn{1}{|c|}{$-$} & 1200 &
\multicolumn{1}{|c}{14} & \multicolumn{1}{|c|}{510} & 4000 &
\multicolumn{1}{|c|}{16000}\\\hline\hline
$\mathcal{I}_{3}(4)$ & no & \multicolumn{1}{|c}{0.30} &
\multicolumn{1}{|c}{0.05} & \multicolumn{1}{|c|}{0.22} & 0.15 &
\multicolumn{1}{|c}{0.03} & \multicolumn{1}{|c|}{0.22} & 1.0 &
\multicolumn{1}{|c|}{8.6}\\\hline
$\mathcal{I}_{3}(5)$ & yes & \multicolumn{1}{|c}{0.72} &
\multicolumn{1}{|c}{0.10} & \multicolumn{1}{|c|}{0.22} & 0.38 &
\multicolumn{1}{|c}{0.07} & \multicolumn{1}{|c|}{0.22} & 110 &
\multicolumn{1}{|c|}{300}\\\hline
$\mathcal{I}_{3}(6)$ & yes & \multicolumn{1}{|c}{1.3} &
\multicolumn{1}{|c}{0.18} & \multicolumn{1}{|c|}{0.22} & 0.83 &
\multicolumn{1}{|c}{0.11} & \multicolumn{1}{|c|}{0.22} & 2500 &
\multicolumn{1}{|c|}{2300}\\\hline\hline
$\mathcal{I}_{4}(6,3)$ & no & \multicolumn{1}{|c}{0.02} &
\multicolumn{1}{|c}{0.01} & \multicolumn{1}{|c|}{0.14} & 0.02 &
\multicolumn{1}{|c}{0.01} & \multicolumn{1}{|c|}{0.14} & 3.1 &
\multicolumn{1}{|c|}{34}\\\hline
$\mathcal{I}_{4}(7,3)$ & no & \multicolumn{1}{|c}{0.04} &
\multicolumn{1}{|c}{0.02} & \multicolumn{1}{|c|}{0.14} & 0.04 &
\multicolumn{1}{|c}{0.01} & \multicolumn{1}{|c|}{0.14} & 1600 &
\multicolumn{1}{|c|}{4000}\\\hline
$\mathcal{I}_{4}(7,4)$ & no & \multicolumn{1}{|c}{0.10} &
\multicolumn{1}{|c}{0.02} & \multicolumn{1}{|c|}{0.14} & 0.10 &
\multicolumn{1}{|c}{0.02} & \multicolumn{1}{|c|}{0.14} & 4300 &
\multicolumn{1}{|c|}{4000}\\\hline
\end{tabular}
\label{tab timings}%
\end{table}

We now turn to the specific examples. In the table, we first consider curves: For  $d\geq1$, let $\mathcal{I}_{1}(d)\in\mathbb{Q}[x_0,\ldots,x_{d}]$ be the ideal generated by
the $2\times 2-$minors of the matrix 
$$
\left(\begin{matrix}
x_0 & \ldots & x_{d-1}\\
x_1 & \ldots & x_d
\end{matrix}\right)
$$
defining a rational
normal curve of degree $d$ in $\mathbb{P}^{d}$. Note that these curves are smooth.

With respect to surfaces, we consider the following examples: The ideals
$\mathcal{I}_{2}(c)\in \mathbb{F}_{103}[x_0,\ldots,x_3,s_1,\ldots,s_{c-1}]$ for $c=2,\ldots,7$ define del Pezzo surfaces of degree $2+c$ of
codimension $c$ in $\mathbb{P}^{2+c}$, as constructed in Example
\ref{ex del pezzo}. The ideals $\mathcal{I}_{3}(m)$ for $m=2,\ldots,6$ correspond to the
unprojection sequence from Example \ref{ex veronese} defined over the rationals. 

Finally the ideals $\mathcal{I}_{4}(d,n)$ define the Stanley-Reisner rings of the cyclic
polytopes of dimension $d$ with $n$ vertices from Example \ref{ex cyclic}, after a random 
linear coordinate change over the rationals with coefficients of bitlength $4$.

\begin{rem}
A parallel implementation of the Jacobian criterion is in the process of
development. The results
obtained so far show that, although the step computing the minors can easily be
parallelized, combining the results of the cores in use into a single ideal 
leads to a significant communication overhead.\footnote{\samepage{The 
parallelization framework for this massive parallel approach is GPI-space, a 
low latency communication library and runtime system for scalable real-time 
parallel applications developed at the Fraunhofer ITWM \cite{GPI}. 
Its use for various computational
settings in algebraic geometry is subject to an ongoing joint project between
ITWM and the {\sc Singular} development team.}} 
Parallelization of Gr\"obner basis
computations is a topic of current research, but still has no general
solution. When computing over the rationals, one approach to parallelization
is to consider modular methods, see \cite{Arnold, FareyPaper}. The conclusion 
is that the Jacobian criterion will benefit much less from parallel 
computations than our algorithm, hence we do not give any timings with regard 
to a parallel Jacobian criterion. 

Note however, that in many applications, for example
when computing the normalization of an affine algebra as in \cite{BDLPSS}, knowledge of
the Jacobian ideal will be required. Thus, even taking the new smoothness test
into account, an efficient parallel approach for the computation
of the Jacobian ideal is of key importance.

\end{rem}

To summarize, we observe that the new algorithm is significantly faster than the Jacobian criterion
in those examples in consideration, where the variety is singular, even when only running on a
single core. For testing smoothness of non-singular varieties, it can be slow if the Groebner basis computations 
involved are complicated, or a very large number of charts has to be considered. Even in these cases, the 
memory requirements are significantly lower than in the Jacobian criterion approach.

The algorithm can benefit to a great extent from parallel and massive parallel
computations, which are of increasing importance considering current and future hardware. As outlined in Section 
\ref{sec hybrid}, this can also be exploited to implement a divide-and-conquer strategy by first doing several descends of ambient dimension using the new approach and then applying the Jacobian criterion in each chart. 
Combining the advantages of both approaches, this hybrid algorithm eventually 
turns out to be the most successful in examples in consideration. 

\vspace{0.1in} \noindent\emph{Acknowledgements}. We would like to thank
Mirko Rahn, Lukas Ristau, Stavros Papadakis, Bernd Schober, and Yue Ren 
for helpful discussions.

\end{document}